\documentclass[12pt,twoside,reqno]{amsart}
\usepackage{pifont}
\usepackage{amsmath}
\usepackage{amsthm}
\usepackage{enumitem}
\usepackage{amsfonts}
\usepackage{amssymb}
\usepackage{latexsym}
\usepackage{amstext}
\usepackage{array}
\usepackage{graphicx}
\usepackage{tikz}
\usepackage{mathrsfs}
\usepackage{diagbox}
\usepackage{threeparttable}
\usepackage{extarrows}

\date{}
\allowdisplaybreaks[4] \footskip=15pt
\renewcommand{\uppercasenonmath}[1]{}

\newtheorem{thm}[subsection]{Theorem\;}
\newtheorem{cor}[subsection]{Corollary\;}
\newtheorem{Def}[subsection]{Definition\;}
\newtheorem{lem}[subsection]{Lemma\;}
\newtheorem{remark}{Remark\;}

\newtheorem{prop}[subsection]{Proposition\;}
\newtheorem{exm}[subsection]{Example\;}

\newcommand{\bthm}{\begin{thm} }
	\newcommand{\ethm}{\end{thm} }
\newcommand{\bpro}{\begin{prop}}
	\newcommand{\epro}{\end{prop}}
\newcommand{\bdf}{\begin{Def}}
	\newcommand{\edf}{\end{Def}}
\newcommand{\bexm}{\begin{exm}}
	\newcommand{\eexm}{\end{exm}}
\newcommand{\blem}{\begin{lem}}
	\newcommand{\elem}{\end{lem}}
\newcommand{\bpf}{\begin{proof}}
	\newcommand{\epf}{\end{proof}}
\newcommand{\bcor}{\begin{cor}}
	\newcommand{\ecor}{\end{cor}}
\newcommand{\ba}{\begin{array}}
	\newcommand{\ea}{\end{array}}
\newcommand{\bea}{\begin{eqnarray}}
	\newcommand{\eea}{\end{eqnarray}}

\newcommand{\brem}{\begin{remark}}
	\newcommand{\erem}{\end{remark}}

\def\deg{\mathrm{deg}}

\graphicspath{dir-list}
\begin{document}
\begin{center}
    {\large  \bf  On the inverse stability of $z^n+c.$ }
  	\footnote {Supported by NSFC (Nos. 12071209, 12231009).}\\
		\vskip 0.8cm
		{\small   Yang Gao and Qingzhong Ji \footnote{Corresponding author.\\ \indent E-mail addresses: 864157905@qq.com (Y. Gao),\; qingzhji@nju.edu.cn (Q. Ji)}}\\
		{\small School of Mathematics, Nanjing University, Nanjing
			210093, P.R.China}
	\end{center}

	{\bf Abstract.} 
     {\small  This paper investigates the inverse stability of the binomials $\phi(z) = z^d + c$ over $K$.
}


\section{\bf  Introduction }\label{01}
Arithmetic dynamics focuses on number-theoretic and algebraic geometric problems arising from iteration. Let \( K \) be a field. A polynomial \( \phi(z) \in K[z] \) is said to be stable if all its iterates are irreducible over \( K \). This concept originated in 1985 with R.W.K. Odoni \cite{Od}, who used the Chebotarev density theorem to provide an asymptotic estimate for the number of prime divisors in the sequence defined by \( a_{n+1} = \phi(a_n) \), with \( a_1 = 2 \) and \( \phi(x) = x^2 - x + 1\in\mathbb{Z}[x] \). The application of the Chebotarev density theorem requires studying the Galois groups of iterated polynomials, which in turn necessitates understanding the reducibility of polynomial iterates.

Subsequently, the concept of stability has been extensively developed by numerous researchers. For instance,  the stability of quadratic polynomials, binomial polynomials, trinomial polynomials, Eisenstein polynomials, and the estimation of the number of stable polynomials over finite fields have all been studied extensively.

\rm{ (See, e.g., \cite{NAli}, \cite{Momc}, \cite{L.Dan},\cite{Doming},\cite{IP},\cite{RL},\cite{Ostafe},\cite{Shanta1},\cite{LW}.)}

In 2017, the concept of polynomial stability was later generalized by R. Jones and L. Alon \cite{Rafe} to the broader notion of eventual stability.
They proposed the ``Everywhere Eventual Stability Conjecture" and demonstrated that it follows from the Dynamical Lehmer Conjecture.

The applications of eventual stability include: (See \cite{Rafe} for details)
\begin{itemize}
    \item The finiteness of \( S \)-relative integral points in backward orbits.  The finiteness of integral points in orbits under a rational function is a dynamical analog of Siegel's theorem on the finiteness of integral points on elliptic curves.
    \item A coarse measure of the image size in the Arboreal Galois representation, which parallels the Galois representation constructed from Tate modules of elliptic curves.
\end{itemize}

Recently, in 2024, K. Cheng \cite{Cheng}  introduced a related concept called inversely stable polynomials. Cheng demonstrated that a polynomial $\phi(z) = z^p + a z + b \in \mathbb{F}_p[z]$ is inversely stable over $\mathbb{F}_p$ if and only if $a = -1$ and $b \neq 0$. 

Moreover, it is straightforward to verify that if $\phi(z) \in K[z]$ is inversely stable over $K$, then $\left(\frac{1}{\phi(z)}, \infty\right)$ is eventually stable over $K$.

In arithmetic dynamics, the study of general phenomena often begins with specific classes of rational maps, such as quadratic polynomials, Lattès maps, Chebyshev polynomials, and binomial polynomials. The iterative behavior of binomial polynomials is an important topic in arithmetic dynamics. For example, it is closely related to problems such as the finiteness of primitive prime divisors in orbits generated by binomial polynomials (see \cite{HK}), height estimates for binomial polynomials (see \cite{PG1}), and the existence of rational periodic points of binomial polynomials (see \cite{PG2}, Theorem 4).

 In this paper, we present several conditions for the binomial polynomial $\phi(z) = z^d + c$ to be inversely stable over the rational number field, function fields, and finite fields, respectively.
In addition, we will give four applications of inverse stability in Section 6. 
\section{\bf  Main result }\label{01}

\begin{Def}
Let $K$ be a field and $\phi(z)\in K[z]$ be a polynomial. Define $\Phi(z) := \frac{1}{\phi(z)} \in K(z).$  For $n \in\mathbb{N}^* $,  let the $n$-th iterate of $\Phi(z)$ be defined as $\Phi^{(n)}(z) = \underbrace{\Phi \circ \Phi \circ \cdots \circ \Phi}_{n \text{ times}}(z).$ We express the \(\Phi^{(n)}(z)\) in its reduced form as
\(
\Phi^{(n)}(z) = \frac{f_{n,\phi}(z)}{g_{n,\phi}(z)},
\)
where \(f_{n,\phi}(z)\) and \(g_{n,\phi}(z)\) are coprime polynomials in \(K[z]\).
 A polynomial $\phi(z) \in K[z]$ is called inversely stable over $K$ if every $g_{n,\phi}(z)$ in the sequence $\{g_{n,\phi}(z)\}_{n=1}^\infty$ is irreducible in $K[z]$.
\end{Def}
\bthm\label{88}
Let \( R \) be a unique factorization domain, and let \( U(R) \) denote the unit group of \( R \). Let \( d \in \mathbb{N}^{*} \) with \( d \geq 2 \), and suppose \( c \notin uR^p \) for all primes \( p \mid d \) and \( u \in U(R)\). Let \( K \) be the fraction field of \( R \). If the polynomial \( \phi(z) = z^d + c \in R[z]\) is irreducible over \( K \), then \( \phi(z)\) is inversely stable over \( K \).
\ethm
\begin{cor}\label{37}
 Let \( d \in \mathbb{N}^{*} \) with \( d \geq 2 .\)
Let \( c \in \mathbb{Z} \) and \( \phi(z) = z^d + c \) be irreducible over \( \mathbb{Q} \). Then \( \phi(z) \) is inversely stable over \( \mathbb{Q} \) if
\begin{itemize}
    \item[(i)] \( d \) is odd, or
    \item[(ii)] \( d \) is even and \( c \) is not a square.
\end{itemize}
\end{cor}

\bthm\label{018}
Let \( K = F(t) \) be the rational function field in one variable over a field \( F \) of characteristic \( 0 \). Let \( d \geq 3 \), and let \( c \in R = F[t] \) with \( c \notin F \). Suppose \( \phi(z) = z^d + c \) is irreducible over \( K \). Then \( \phi(z) \) is inversely stable over \( K \).
\ethm
\bthm\rm{(\cite{RafeB}, Proposition 2.3.)} Let \( K \) be a finite field of characteristic not equal to two. A quadratic polynomial \( f(z)=az^2+bz+c \in K[z] \) is stable if and only if the set 
\[
\{-f(\gamma)\} \cup \left\{ f^i(\gamma) : i = 2, 3, \ldots \right\}
\]
contains no squares, where $\gamma=-\frac{b}{2a}.$ 
\ethm
In this paper, we establish the following similar result.
\bdf \rm{(\cite{SH}, Definition 5.1)}  For $m \mid q-1, \;\alpha \in \mathbb{F}_q^*$ is $m$-free if the equality $\alpha=\beta^d$ with $\beta \in \mathbb{F}_q$, for any divisor $d$ of $m$, implies $d=1$.
\edf
\bthm\label{20}
Let \( K = \mathbb{F}_q \). Let \( \phi(z) = z^d + c \in K[z] \), where \( d \geq 2 \). Suppose that \( \phi(z) \) is irreducible over \( K \).

Define the sequence:
\[
x_1 = c, \quad x_2 = (-1)^d(c^{d+1} + 1), \quad x_{n+2} = (-1)^d c x_{n+1}^d + x_n^{d^2}, \quad n \in \mathbb{N}^*.
\]

Then \( \phi \) is inversely stable over \( K \) if and only if \( \frac{x_{n+1}}{x_n} \) is \( \operatorname{rad}(d) \)-free for every \( n \in \mathbb{N}^* \), where $
\operatorname{rad}(d)=\prod\limits_{p \mid d, \;p \text { is a prime }} p
$.
\ethm
\begin{cor}\label{21}
Let \( n \in \mathbb{N}^* , n\geq 3\) , $d=2^n$ and suppose \( p = 2^{n+1} + 1 \) is a prime. Define 
\[
\phi(z) = z^{d} + c \in \mathbb{F}_{p}[z].
\]
Then there are at least 
\(
\lfloor \frac{1}{8} \left( \sqrt{p} - 1 \right)^2\rfloor
\)
distinct values of \( c \in \mathbb{F}_{p} \) such that \( \phi(z) \) is inversely stable over \( \mathbb{F}_{p} \), where $\lfloor x \rfloor$ denotes the greatest integer less than or equal to $x.$ 
\end{cor}

This paper is organized as follows: In $\S 3,$ we shall give the proof of Theorem {\rm\ref{88}}. In $\S 4,$ 
we shall give the proof of Theorem {\rm\ref{018}}. In $\S 5,$ we shall give the proof of Theorem {\rm\ref{20}} and Corollary {\rm\ref{21}}.
In $\S 6$, we conclude the paper by summarizing the significance of our work and presenting its potential applications.

\section{\centering \textbf{Proof of Theorem {\rm\ref{88}}}}
Let $K$ be a field.
		A rational function $\varphi(z)=\frac{f(z)}{g(z)} \in K(z)$ is a quotient of polynomials $f(z), g(z)\in K[z]$ with no common factors.
		 The degree of $\varphi$ is
		$
		\operatorname{deg} \varphi=\max \{\operatorname{deg} f, \operatorname{deg} g\}.
		$
		The rational function $\varphi$ of degree $d$ induces a rational map (morphism) of the projective space $\mathbb{P}^1(\overline{K}),$
		$$
		\varphi: \mathbb{P}^1(\overline{K}) \longrightarrow \mathbb{P}^1(\overline{K}),\;\;
		\varphi([X:Y]) = \left[Y^d f(X / Y): Y^d g(X / Y)\right].$$ A point $P \in \mathbb{P}^1(\overline{K})$ is said to be periodic under $\varphi$ if $\varphi^{(n)}(P) = P$ for some $n \geq 1$.

  We write \([1:0] \in \mathbb{P}^1(K)\) as \(\infty\).
\blem\label{002}
Let $d \in \mathbb{N}^*$, and let $K$ be a field such that $\operatorname{char}(K) = 0$ or $\operatorname{char}(K) > 0$ with $\operatorname{char}(K)$ prime to $d$. Consider the polynomial $\phi(z) = z^d + c \in K[z], c\neq 0$ and define the rational function $\Phi(z) := \frac{1}{\phi(z)} \in K(z)$. For each $n \in \mathbb{N}^*$, denote the $n$-th iterate of $\Phi$ by $\Phi^{(n)}(z) = \frac{f_{n,\phi}(z)}{g_{n,\phi}(z)}$, where $f_{n,\phi}(z)$ and $g_{n,\phi}(z) \in K[z]$ are coprime polynomials. If $\infty$ is not periodic under $\Phi(z)$, then for any $n \in \mathbb{N}^*$, the degree of $g_{n,\phi}(z)$ is $d^n$.
\elem
\bpf
Note that the map \( \Phi^{(n)} : \mathbb{P}^1(\overline{K}) \to \mathbb{P}^1(\overline{K}) \) is given by
\[
\Phi^{(n)}([X : Y]) = \left[Y^e f_{n, \phi}\bigg( \frac{X}{Y} \bigg) : Y^e g_{n, \phi}\bigg( \frac{X}{Y} \bigg)\right],
\]
where \( e = \deg \Phi^{(n)}(z) \).

It follows that $\Phi^{(n)}([\alpha : 1]) = \left[f_{n, \phi}(\alpha) : g_{n, \phi}(\alpha)\right]$, and hence
\[
\Phi^{(n)}([\alpha : 1]) = \infty \quad \text{if and only if} \quad g_{n, \phi}(\alpha) = 0.
\]
By assumption we have \( \infty \notin (\Phi^{(n)})^{-1}(\infty) \) for any \( n \in \mathbb{N}^* \). Thus, for \( n \geq 1 \),
\[
(\Phi^{(n)})^{-1}(\infty) = \{ [\alpha : 1] \in \mathbb{P}^1(\overline{K}) \mid \Phi^{(n)}([\alpha : 1]) = \infty \}.
\]
Thus,
\[
(\Phi^{(n)})^{-1}(\infty) = \{ [\alpha : 1] \in \mathbb{P}^1(\overline{K}) \mid g_{n, \phi}(\alpha) = 0 \}.
\]
Next, we prove that \( \# (\Phi^{(n)})^{-1}(\infty) = d^n \).

Since
\(
\Phi([X : Y]) = [Y^d : X^d + c Y^d],
\)
we have \( \Phi(\infty) = [0 : 1] \) and \( \Phi([0 : 1]) = [1 : c] \).
Since \( \infty \notin (\Phi^{(n)})^{-1}(\infty) \) for all \( n \), neither \( [0 : 1] \) nor \( [1 : c] \) belong to \( (\Phi^{(n)})^{-1}(\infty) \). 

For any $P = [1 : t] \in \mathbb{P}^1(\overline{K})$, we have  
\[
\Phi([X : Y]) = P \quad \text{if and only if} \quad \left[Y^d : X^d + c Y^d\right] = [1 : t],
\]  
which simplifies to $\left(\frac{X}{Y}\right)^d + c - t = 0$. Therefore, if $t \neq c$, it follows that $\#\Phi^{-1}(P) = d$.

Hence,  
\[
\left|\left(\Phi^{(1)}\right)^{-1}(\infty)\right| = d \quad \text{and} \quad \left|\left(\Phi^{(i+1)}\right)^{-1}(\infty)\right| = d \left|\left(\Phi^{(i)}\right)^{-1}(\infty)\right|
\]  

for any \(i\in\mathbb{N}^*\). It follows that  
\(
\left|\left(\Phi^n\right)^{-1}(\infty)\right| = d^n.
\)
Thus, \(g_{n,\phi}\) has \(d^n\) distinct roots in \(\overline{K}\). Combining this result with \(\deg(g_{n,\phi}) \leq d^n\), we conclude that \(\deg(g_{n,\phi}) = d^n\) for any $n\in\mathbb{N}^*$.
\qed\epf
\blem\label{003}
Let \( F \) be a field, and let \( f(z) = z^d + m \in F[z] \) be an irreducible polynomial.  
Denote by \( \overline{F} \) the algebraic closure of \( F \), and let \( \gamma \in \overline{F} \) be a root of \( f(z) \).  
Let \( a, b, e, t \in F \) with \( ae \neq 0 \).  
We denote by \( N_{F(\gamma)/F} \) the norm map associated with the field extension \( F(\gamma)/F \).
 Then 

\[
N_{F(\gamma)/F}\left(\frac{a\gamma+b}{e\gamma+t}\right) = \frac{b^d + (-1)^d m a^d}{t^d + (-1)^d m e^d}.
\]
\elem
\bpf
The conjugates of \( \frac{a\gamma + b}{e\gamma + t} \) are \( \frac{a\gamma_i + b}{e\gamma_i + t} \) for \( i = 1, 2, \dots, d \), where \( \gamma_1, \gamma_2, \dots, \gamma_d \in\overline{F}\) are the roots of \( f(z) = z^d + m \). The norm is
\[
N_{F(\gamma)/F} \left( \frac{a\gamma + b}{e\gamma + t} \right) = \prod_{i=1}^d \frac{a\gamma_i + b}{e\gamma_i + t}.
\]
Using the fact that $f(z)=\prod\limits_{i=1}^d(z-\gamma_i),$ we have
\[
\prod_{i=1}^d (a\gamma_i + b) = a^d (-1)^d m + b^d, \quad \prod_{i=1}^d (e\gamma_i + t) = e^d (-1)^d m + t^d.
\]
Therefore, 
\[
N_{F(\gamma)/F} \left( \frac{a\gamma + b}{e\gamma + t} \right) = \frac{b^d + (-1)^d m a^d}{t^d + (-1)^d m e^d}.
\]
This completes the proof.
\qed
\epf
\blem\label{110}
Let \(d \in \mathbb{N}^*\) with \(d \geq 2\), and let \(c \in R \setminus \{0\}\), where \(R\) is a unique factorization domain, and \(c\) is a non-unit element of \(R\).
 Define a sequence of matrices \( \{A_j\}_{j \geqslant 1} \) in \( M_{2 \times 2}(R) \) by the following relations:

\[
A_1 = \begin{bmatrix} 
x_1 & y_1 \\
z_1 & w_1 
\end{bmatrix} 
= \begin{bmatrix}
c & -1 \\
1 & 0 
\end{bmatrix},
\]
and for \( j \geq 1 \),
\[
A_{j+1} = \begin{bmatrix} 
x_{j+1} & y_{j+1} \\
z_{j+1} & w_{j+1} 
\end{bmatrix} 
= \begin{bmatrix}
(-1)^d c x_j^d + y_j^d & (-1)^{d+1} x_j^d \\
(-1)^d c z_j^d + w_j^d & (-1)^{d+1} z_j^d
\end{bmatrix}.
\]

Then the following statements hold:

\begin{enumerate}[label=\textup{(\roman*)}]
    \item For all \( n \geq 1 \),
    \[
   x_{n+2} = (-1)^d c x_{n+1}^d + x_n^{d^2}, \quad \gcd(x_{n+1}, x_n) = 1,\quad z_{n+1} = (-1)^d x_n.
    \]
    
    \item For all \( n \geq 1 \),
    \[ c\mid x_{2n-1},\;c\mid x_{2n}-(-1)^{d},\;
    \gcd\left( \frac{x_{2n-1}}{c}, c \right) = 1.
    \]
    
    \item If \( c \notin uR^p \) for all primes \( p \mid d \) and units \( u \in R \), then for all \( n \geq 1 \),
    \[
    x_{2n-1} \notin uR^p \quad \text{for all} \; p \mid d \text{ and units} \; u \in R.
    \]
   \item  \( \infty \) is not periodic under \( \Phi(z)=\frac{1}{z^d+c} \) if and only if 
    \(x_n \neq 0\) for all \( n \geq 1 .\)
    \item Assume that \( \infty \) is not periodic under the map \( \Phi(z) = \frac{1}{z^d + c} \). If \( c \notin u R^p \) for all primes \( p \mid d \) and units \( u \in R \), then for all \( n \geq 1 \), we have
\[
\frac{x_{n+1}}{x_n} \notin \pm K^p \quad \text{for all primes} \; p \mid d, 
\]
where $K$ is the fraction field of $R.$

\end{enumerate}
\elem
\bpf  
 (i) and (ii) are trivial from definition. 
 
 (iii).
Assume that \( x_{2n-1} = u_1 r_1^{p_1} \) for some unit \( u_1 \in R \), \( r_1 \in R \) and prime \( p_1 \mid d \). Then, we have  
\[
c \cdot \frac{x_{2n-1}}{cu_1} = r_1^{p_1}.
\]  
By (ii), we know that \( c \) and \( \frac{x_{2n-1}}{cu_1} \) are coprime elements in \( R \).  
It follows that \( c \) can be written as \( c = u_2 r_2^{p_1} \), where \( u_2 \in R \) is a unit and \( r_2 \in R \).  
This contradicts the assumption.

(iv). Define the sequences \( \{a_n\}_{n\in\mathbb{N}^*} \) and \( \{b_n\}_{n\in\mathbb{N}^*} \) in \( R \) as follows:  
\[
a_1 = 0, \quad b_1 = 1, \quad a_{n+1} = b_n^d, \quad b_{n+1} = a_n^d + c b_n^d.
\]  
Then, \( \Phi^{(n)}(\infty) = [a_n : b_n] \). 

Now, observe that  
\[
b_{n+2} = c b_{n+1}^d + b_n^{d^2}, b_1 = 1, b_2 = c.
\]  
It is obvious that  
\(
x_{2n-1} = b_{2n}, x_{2n} = (-1)^d b_{2n+1}.
\)
Hence, \( \infty \) is not periodic under \( \Phi \) if and only if \( x_n \neq 0 \) for all \( n \geq 1\).  

(v).
Assume \( \frac{x_{n+1}}{x_n} \in \pm K^p \) for some prime \( p \mid d \). Since \( \gcd(x_{n+1}, x_n)= 1 \), it follows that there exist units \( u_3, u_4 \in R \) such that  
\[
x_{n+1} \in u_3 R^p \quad \text{and} \quad x_n \in u_4 R^p.
\]  
 By (iii), this situation is impossible.
 \qed  
\epf

\begin{lem} \label{001}\rm{(\cite{GK}, Theorem 8.1.6.)}
Let \( K \) be a field, \( d \geq 2 \) an integer, and \( a \in K \). The polynomial \( X^d + a \) is irreducible over \( K \) if and only if \( a \notin -K^p \) for all primes \( p \) dividing \( d \), and \( a \notin 4K^4 \) whenever \( 4 \mid d \).
\end{lem}

\textbf{Proof of Theorem {\rm\ref{88}}}
Let the sequence $\{x_n\}_{n\in\mathbb{N}^*}$ be defined in Lemma~{\rm{\ref{110}}}.
 Claim 1:  \( \infty \) is not periodic under \( \Phi \).

If \( -c \in K^d \), then \( z^d + c = z^d - (-c) \) is reducible over \( K \). Therefore, \( -c \notin K^d \).

 We verify the initial terms of the sequence $\{x_n\}_{n\in\mathbb{N}^*}$,
\( x_1 = c \neq 0 \) and \( x_2 = (-1)^d (c^{d+1} + 1) \neq 0 \); otherwise, we would have \( -c = \left(\frac{1}{c}\right)^d \in K^d\).

 Assume \( x_n \neq 0 \) and \( x_{n+1} \neq 0 \). If \( x_{n+2} = 0 \), then \( -c = \left(-\frac{x_n^d}{x_{n+1}}\right)^d \in K^d\). This contradicts the assumption  \( -c \notin K^d \). Hence, \( x_n \neq 0 \) for all \( n \in \mathbb{N}^* \). 

By Lemma~{\rm{\ref{110}}}~{\rm{(iv)}}, \( \infty \) is not periodic under \( \Phi \). This completes the proof of the Claim 1.

Let \(\{Q_i\}_{i \geq 1}\) be a sequence in \(\mathbb{P}^1(\overline{K})\) such that \(\Phi(Q_1) = \infty\) and \(\Phi(Q_{i+1}) = Q_i\) for all \(i \geq 1\). Since \(\infty\) is not periodic under \(\Phi\), and \(\Phi(\infty) = [0:1]\), we can express each \(Q_i\) as \(Q_i = [\beta_i : 1]\), where \(\beta_i \in \overline{K}\) and \(\beta_i \neq 0\) for all \(i \in \mathbb{N}^{*}\).

Thus, we have
\[
\phi(\beta_1) = 0, \phi(\beta_{i+1}) = \frac{1}{\beta_i}.
\]
It is clear that \(\beta_n\) is a root of the polynomial \(g_{n,\phi}(z)\).

 Claim 2:  \( z^d + c - \frac{1}{\beta_n} \) is irreducible over \( K(\beta_n) \) for every $n\geq 1$. 
 
 We shall prove this claim by induction on $n.$

By Lemma~{\rm{\ref{003}}}, we have  
\[
N_{K(\beta_1)/K}\left(\frac{c\beta_1 - 1}{\beta_1}\right) = \frac{c^{d+1} + 1}{c} = (-1)^d \frac{x_2}{x_1},
\]
where \( x_1 \) and \( x_2 \) are as defined in Lemma~{\rm{\ref{110}}}.

By Lemma~{\rm{\ref{110}}}~(v), we deduce that \( (-1)^d \frac{x_2}{x_1} \not\in \pm K^p \) for all primes \( p \mid d \), and hence 
\[
\frac{c\beta_1 - 1}{\beta_1} \not\in -K(\beta_1)^p \quad \text{for all primes } p \mid d.
\]
Obviously, if \( 4 \mid d\) and \( \frac{c\beta_1 - 1}{\beta_1} \in 4K(\beta_1)^4 \), then \( \frac{c\beta_1 - 1}{\beta_1} \in K(\beta_1)^2 \) and \( (-1)^d \frac{x_2}{x_1} \in  K^2 \), which also contradicts Lemma~{\rm{\ref{110}}}~{\rm{(v)}}. So \( \frac{c\beta_1 - 1}{\beta_1} \notin 4K(\beta_1)^4 \), when $4\mid d.$
 By Lemma~{\rm{\ref{001}}}, the Claim 2 holds for $n=1$. 

Therefore, \( [K(\beta_2) : K(\beta_1)] = d \).

Assume that $[K(\beta_{i}):K(\beta_{i-1})]=d$ for each \( 2 \leq i \leq n .\)
We will prove that $[K(\beta_{n+1}):K(\beta_{n})]=d.$
This means that we will prove that
\( z^d + c - \frac{1}{\beta_{n}} \) is irreducible over the field \( K(\beta_n) .\)

Based on the inductive hypothesis, we know that  \( z^d + c - \frac{1}{\beta_j} \) is irreducible over the field \( K(\beta_j) \) for each \( j \) with \( 1 \leq j \leq n-1 \).

Given a fixed $j $ with \( 1 \leq j \leq n-1 \). Set $F=K(\beta_j), \gamma=\beta_{j+1}, f(z)=z^d+c-\frac{1}{\beta_j}.$

Then  \( K(\beta_{j+1}) = K(\beta_j)(\beta_{j+1}) =F(\beta_{j+1}).\)  Therefore, we obtain 
\begin{align*}
    N_{K(\beta_{j+1})/K(\beta_j)}\left(\frac{x_l \beta_{j+1} + y_l}{z_l \beta_{j+1} + w_l}\right) 
    &= N_{F(\gamma)/F}\left(\frac{x_l \gamma + y_l}{z_l \gamma + w_l}\right) 
    &\overset{\text{Lemma~{\rm{\ref{003}}}}}{=}\frac{y_l^d + (-1)^d(c - \frac{1}{\beta_j})x_l^d}{w_l^d + (-1)^d(c - \frac{1}{\beta_j})z_l^d}, 
    \end{align*}
    and therefore  $$N_{K(\beta_{j+1})/K(\beta_j)}\left(\frac{x_l \beta_{j+1} + y_l}{z_l \beta_{j+1} + w_l}\right) = \frac{x_{l+1} \beta_j + y_{l+1}}{z_{l+1} \beta_j + w_{l+1}},$$
for any \( 1 \leq j \leq n-1 \) and \( l \in \mathbb{N}^* \),
where \( x_l \), \( y_l \), \( z_l \), and \( w_l \) are the terms of index \( l \) in the sequence defined in Lemma~{\rm{\ref{110}}}. (Note that we have previously proved that \( x_n = (-1)^d z_{n+1} \neq 0 \) for any $n\in\mathbb{N}^{*}$.)

This implies that 
\[
N_{K(\beta_{n})/K(\beta_{n-1})}\left(\frac{x_1\beta_{n}+y_1}{z_1\beta_{n}+w_1}\right) = \frac{x_2\beta_{n-1}+y_2}{z_2\beta_{n-1}+w_2}
\]
\[
N_{K(\beta_{n-1})/K(\beta_{n-2})}\left(\frac{x_2\beta_{n-1}+y_2}{z_2\beta_{n-1}+w_2}\right) = \frac{x_3\beta_{n-2}+y_3}{z_3\beta_{n-2}+w_3}
\]
\[
\vdots
\]
\[
N_{K(\beta_{2})/K(\beta_{1})}\left(\frac{x_{n-1}\beta_{2}+y_{n-1}}{z_{n-1}\beta_{2}+w_{n-1}}\right) = \frac{x_{n}\beta_{1}+y_{n}}{z_{n}\beta_{1}+w_{n}}
\]
By  Lemma~{\rm{\ref{003}}}  and the fact $\phi(z)=z^d+c$ is irreducible over $K,$ 
we obtain 
\[
N_{K(\beta_{1})/K}\left(\frac{x_n\beta_{1}+y_n}{z_n\beta_{1}+w_n}\right) = \frac{x_{n+1}}{z_{n+1}} = (-1)^d \frac{x_{n+1}}{x_{n}}
.\]
So, 
\[
N_{K(\beta_{n})/K}\left(\frac{x_1\beta_{n}+y_1}{z_1\beta_{n}+w_1}\right) = (-1)^d \frac{x_{n+1}}{x_{n}}
.\]
Finally, by Lemma~{\rm{\ref{110}}}~{\rm{(v)}}, we conclude that
\[
c - \frac{1}{\beta_n}=\frac{x_1\beta_{n}+y_1}{z_1\beta_{n}+w_1} \notin  -K(\beta_n)^p \quad \text{for all primes} \quad p \mid d,
\]
and
\[
c - \frac{1}{\beta_n} \notin  4K(\beta_n)^4 \quad \text{whenever} \quad 4 \mid d.
\]
From Lemma~{\rm{\ref{001}}}, we deduce that
\( z^d + c - \frac{1}{\beta_{n}} \) is irreducible over the field \( K(\beta_n) .\)
This completes the proof of the Claim 2.

Therefore, \( [K(\beta_n) : K] = d^n \) for any \( n \in \mathbb{N}^* \). Since \( \beta_n \) is a root of \( g_{n, \phi} \) and, by Lemma~{\rm{\ref{002}}}, we have \( \operatorname{deg}(g_{n, \phi}) = d^n \), it follows that \( g_{n, \phi} \) is irreducible over \( K \) for all \( n \in \mathbb{N}^* \).

\qed

\begin{cor}
Let $c \in \mathbb{Z}$ and let \( \phi(z) = z^d + c \) be irreducible over \( \mathbb{Z} \). Then, \( \phi(z) \) is inversely stable over \( \mathbb{Q} \) if:
\begin{itemize}
    \item[\rm{(i)}] \( d \) is odd, or
    \item[\rm{(ii)}] \( d \) is even and \( c \) is not a perfect square.
\end{itemize}
\end{cor}

\bpf
Let \( R = \mathbb{Z} \). Since \( \phi(z) = z^d + c \in R[z] \) is irreducible, it follows that \( c \notin -R^p \) for any prime \( p \) dividing \( d \).  
When \( p \) is odd, it is clear that \( -R^p = R^p \). Note that \( U(R) = \{\pm 1\} \).  
Applying Theorem~\ref{88}, the proof is complete.  
\qed
\epf

\section{\centering \textbf{Proof of Theorem {\rm\ref{018}}}}
\blem \label{150}\rm{(\cite{SL}, Theorem 7.1)}
Let \( K \) be a field with characteristic \( 0 \), and let \( \overline{K} \) be its algebraic closure. For a polynomial \( f(t) \in K[t] \), define \( n_0(f) \) to be the number of distinct roots of \( f \) in \( \overline{K} \).
Let \( a(t), b(t), c(t) \in K[t] \) be polynomials that are relatively prime, such that \( a(t) + b(t) = c(t) \)  and  not all of them have vanishing derivative. Then, we have the inequality
\[
\max \{ \deg(a), \deg(b), \deg(c) \} \leq n_0(a(t) b(t) c(t)) - 1.
\]
\elem
\textbf{Proof of Theorem {\rm\ref{018}}}
\begin{proof}
Define the sequence \( \{x_n\}_{n \in \mathbb{N}^*} \) by:
\[
x_1 = c, \quad x_2 = (-1)^d(c^{d+1}+1), \quad x_{n+2} = (-1)^d c x_{n+1}^d + x_n^{d^2}
\tag{1}
\]
for all \( n \in \mathbb{N}^* \). This sequence is consistent with the sequence \( \{x_n\}_{n \in \mathbb{N}^*} \) described in Lemma~\ref{110}. 

For any \( n \in \mathbb{N}^* \), we have the following degree formula:
\[
\operatorname{deg}(x_n) = \frac{d^n - 1}{d - 1} \operatorname{deg}(c),
\tag{2}
\]
where \( \operatorname{deg}(c) \) denotes the degree of \( c(t) \) viewed as a polynomial in \( t \).

We claim that \( x_{2n} \notin u R^p \) for any unit \( u \in R \) and any prime \( p \) dividing \( d \), for all \( n \geq 1 \). Suppose that there exist \( u \in R^* \), \( k \in \mathbb{N}^* \), \( z_k \in R \), and a prime \( p \mid d \) such that 
\[
x_{2k+2} = u z_k^p.
\tag{3}
\]
We then define
\[
g_k = \frac{1}{u} (-1)^d c x_{2k+1}^d \quad \text{and} \quad h_k = \frac{1}{u} x_{2k}^{d^2}.
\tag{4}
\]
From (1), (2), (3) and (4), we have \( \operatorname{deg}(g_k) \geq 1 \), \( \operatorname{deg}(h_k) \geq 1 \), \( \operatorname{deg}(z_k) \geq 1 \), and 
\[
g_k + h_k = z_k^p.
\tag{5}
\]
By Lemma~\ref{110} (i) and (ii), it follows that 
\[
\gcd(x_{2k}, c) = 1 \quad \text{and} \quad \gcd(x_{2k}, x_{2k+1}) = 1.
\tag{6}
\]
In light of (4), (5), and (6), we obtain \( g_k \), \( h_k \), and \( z_k \) are pairwise coprime. Applying Lemma~\ref{150} on (5), we obtain the following inequality:
\[
\begin{aligned}
\operatorname{deg}(g_k) + \operatorname{deg}(h_k) + \operatorname{deg}(z_k^p) 
&\leq 3 \left( n_0(g_k h_k z_k^p) - 1 \right) \\
&= 3 \left( n_0(g_k) + n_0(h_k) + n_0(z_k) - 1 \right).   
\end{aligned}
\tag{7}
\]

By Lemma~\ref{110}\;(ii), we have \( \gcd\left(c, \frac{x_{2k+1}}{c}\right) = 1 \). Thus, 
\[
\begin{aligned}
n_0(g_k) &= n_0(c^{d+1}) + n_0\left( \left(\frac{x_{2k+1}}{c}\right)^d \right) \\
         &= n_0(c) + n_0\left( \frac{x_{2k+1}}{c} \right) \leq \operatorname{deg}(c) + \operatorname{deg}\left( \frac{x_{2k+1}}{c} \right) \\
         &= \operatorname{deg}(x_{2k+1}),
\end{aligned}
\tag{8}
\]
Additionally, we have
\[
n_0(h_k) = n_0(x_{2k}) = \operatorname{deg}(x_{2k}),\;\text{and}\; 
 n_0(z_k) \leq \operatorname{deg}(z_k) .
\tag{9}
\]
Combining (7), (8), and (9) yields 
\[
(d - 3) \operatorname{deg}(x_{2k+1}) + (d^2 - 3) \operatorname{deg}(x_{2k})
+ (p - 3) \operatorname{deg}(z_k) + 4 \leq 0. 
\tag{10}
\]

If \( d \) is odd, this inequality leads to a contradiction, as \( p \mid d \). Therefore, we need to consider the case where \( d \) is even and \( p = 2 \).

Applying Lemma~\ref{150} again yields the inequality:
\[
p \operatorname{deg}(z_k)= \operatorname{deg}(z_k^p)\leq \operatorname{deg}(x_{2k+1}) + \operatorname{deg}(x_{2k}) + \operatorname{deg}(z_k) - 1. 
\tag{11}
\]

Combining inequalities (10) and (11) with \( p = 2 \), we obtain:
\[
(d - 4) \operatorname{deg}(x_{2k+1}) + (d^2 - 4) \operatorname{deg}(x_{2k}) + 5 \leq 0,
\]
which is a contradiction.
Thus, the claim is proved.
From the proof of Theorem~\ref{88}, it is evident that \( \phi(z) \) is inversely stable over \( K \).\qed
\end{proof}

\section{\centering\textbf{Proof of Theorem {\rm\ref{20}} and Corollary {\rm\ref{21}}}}
\blem\label{009}\rm{(\cite{ADL}, Proposition 2.3)}
Let \( t \geq 2 \) be an integer and \( b \in \mathbb{F}_q \). Then the binomial \( x^t - b \) is irreducible in \( \mathbb{F}_q[x] \) if and only if the following conditions are satisfied:
\begin{enumerate}[label=\textup{(\roman*)}]
    \item \( \operatorname{rad}(t) \mid q - 1 \);
    \item \( b \) is \( \operatorname{rad}(t) \)-free;
    \item \( q \equiv 1 \pmod{4} \) if \( t \equiv 0 \pmod{4} \).
\end{enumerate}
\elem
\blem\label{010}\rm{(\cite{ADL}, Corollary 2.8)}
 Let $t\geq 2$ be an integer such that $\operatorname{rad}(t) \mid q-1$. Then an element $\alpha \in \mathbb{F}_{q^n}$ is $\operatorname{rad}(t)$-free if and only if $N_{q^n / q}(\alpha)$ is $\operatorname{rad}(t)$-free in $\mathbb{F}_q$.
\elem
\textbf{Proof of Theorem {\rm\ref{20}}}

Obviously, \( g_{n, \phi} \) is irreducible if and only if \( [K(\beta_n) : K] = d^n \), where $\{\beta_i\}_{i\in\mathbb{N}^*}$ are described in the proof of Theorem~{\rm{\ref{88}}}. It is easy to observe that \( [K(\beta_i) : K(\beta_{i-1})] \leq d \) for all \( i \geq 2 \), and \( [K(\beta_1) : K] \leq d \).  Therefore, \( g_{n, \phi} \) is irreducible if and only if \( [K(\beta_{i+1}) : K(\beta_{i})] = d \) for all \( 1 \leq i \leq n-1 \), and \( [K(\beta_1) : K] = d \). The latter is equivalent to that \( z^d + c - \frac{1}{\beta_i} \) is irreducible over \( K(\beta_i) \) for all \( 1 \leq i \leq n-1 \), and \( z^d + c \) is irreducible over \( K \).

By Lemma~{\rm{\ref{009}}} and Lemma~{\rm{\ref{010}}}, \( g_{n, \phi} \) is irreducible if and only if the following conditions hold:
\begin{enumerate}[label=\textup{(\roman*)}]
    \item \( \operatorname{rad}(d) \mid q - 1 \);
    \item Both \( -c \) and \( \frac{x_{i+1}}{x_i} \) are \( \operatorname{rad}(d) \)-free  for all \( 1 \leq i \leq n-1 \);
    \item If \( 4 \mid d \), then \( q \equiv 1 \mod 4 \).
\end{enumerate}
Certain details are omitted here for brevity; however, these details are fully addressed in the proof of Theorem~\ref{88}.
\qed
\blem \rm{(\cite{JHS}, Application 1.3, page 139)}\label{50}
 Let $\mathbb{F}_q$ be a finite field with $q$ odd.  Let
$$
f(x)=a x^3+b x^2+c x+d \in \mathbb{F}_q[x]
$$
be a cubic polynomial with distinct roots in $\overline{\mathbb{F}}_q$, and let
$$
\chi: \mathbb{F}_q^* \longrightarrow\{ \pm 1\}
$$
be the unique nontrivial character of order 2, i.e., $\chi(t)=1$ if and only if $t$ is a square in $\mathbb{F}_q^*$. Extend $\chi$ to $\mathbb{F}_q$ by setting $\chi(0)=0$. Then
$$
\left|\sum_{x \in \mathbb{F}_q} \chi(f(x))\right| \leq 2 \sqrt{q} .
$$
\elem
\blem\label{60} \label{118}\rm{(\cite{RH}, Theorem 5.48)}
Assume that $\mathbb{F}_q$ is a finite field of odd characteristic, and assume that $\chi$ is the quadratic character on $\mathbb{F}_q$. Then, if $f(x)=a x^2+b x+c$ is a quadratic polynomial in $\mathbb{F}_q[x]$, the following identity holds :

\[
\sum_{x \in \mathbb{F}_q} \chi\left(a x^2+b x+c\right) = 
\begin{cases}
-\chi(a) & \text{if } b^2 - 4ac \neq 0, \\
\chi(a)(q-1) & \text{if } b^2 - 4ac = 0.
\end{cases}
\]
\elem
\textbf{Proof of  Corollary {\rm\ref{21}}}

Firstly, we claim that if \( \left(\frac{c-1}{p}\right) = 1 \) and \( \left(\frac{c}{p}\right) = \left(\frac{c+1}{p}\right) = -1 \), where \( \left(\frac{\cdot}{p}\right) \) denotes the Legendre symbol, then \( \phi(z) \) is inversely stable over \( \mathbb{F}_p \).
On the one hand, by Lemma~{\rm{\ref{009}}} and the fact that 
\[
\left( \frac{-c}{p} \right) = \left( \frac{-1}{p} \right) \left( \frac{c}{p} \right) = -1, \quad \text{i.e., } c \text{ is 2-free},
\]
we obtain that \( z^d + c \) is irreducible over \( K = \mathbb{F}_p \).  
This implies that \( -c \notin K^d \).  
From the recurrence relation of \( x_n \), it is straightforward to see that \( x_n \neq 0 \) for any \( n \geq 0 \).  

On the other hand,
note that \( p-1 \mid d^2 \). Hence, by Fermat's Little Theorem, for any integer \( a \) with \( p \nmid a \), we have \( a^{d^2} \equiv 1 \pmod{p} \).

Thus, from \( x_n \neq 0 \), we conclude that \( x_n^{d^2} = 1 \).  
From the recurrence relation, we obtain
\[
x_1 = c, \quad x_2 = (-1)^d (c^{d+1} + 1), \quad x_{n+2} = (-1)^d c x_{n+1}^d + 1, \quad n \in \mathbb{N}^*.
\]
 From Euler's Criterion, for any integer \( a \) and an odd prime \( p \), we have
\[
\left( \frac{a}{p} \right) \equiv a^{\frac{p-1}{2}} \pmod{p}.
\]

Hence, 
\[
(c-1)^d = 1, \quad c^d = -1, \quad (c+1)^d = -1.
\]
So we have \( x_2 = 1 - c \), \( x_3 = c + 1 \), \( x_4 = 1 - c \), and \( x_5 = c + 1 \).  
Thus, the sequence \( \left\{ x_n \right\}_{n \geqslant 1} \) follows the pattern:
\[
\left\{ x_n \right\}_{n \geqslant 1}: c, 1 - c, c + 1, 1 - c, c + 1, \ldots
\]
The sequence \( \left\{ \frac{x_{n+1}}{x_n} \right\}_{n \geqslant 1} \) is given by:
\[
\left\{ \frac{x_{n+1}}{x_n} \right\}_{n \geqslant 1}: \frac{1-c}{c}, \frac{1+c}{1-c}, \frac{1-c}{1+c}, \frac{1+c}{1-c}, \frac{1-c}{1+c}, \ldots
\]
From \( \left( \frac{c-1}{p} \right) = 1 \) and \( \left( \frac{c}{p} \right) = \left( \frac{c+1}{p} \right) = -1 \), we know that both \( \frac{1-c}{c} \), \( \frac{1+c}{1-c} \), and \( \frac{1-c}{1+c} \) are 2-free, and that \( -c \) is also 2-free.

By Theorem~{\rm \ref{20}} and Lemma~{\rm \ref{009}}, we conclude that \( \phi(z) \) is inversely stable over \( \mathbb{F}_p \).
This completes the proof of the claim.

To prove that there are at least $\frac{1}{8}(\sqrt{p}-1)^2$ distinct values of $c \in \mathbb{F}_p$ such that $\left(\frac{c-1}{p}\right)=1$ and $\left(\frac{c}{p}\right)=\left(\frac{c+1}{p}\right)=-1$, we analyze the function $t(x)$ defined as:

$$
t(x)=\frac{(\chi(x-1)+1)}{2} \cdot \frac{(1-\chi(x))}{2} \cdot \frac{(1-\chi(x+1))}{2}
$$

where $\chi(x)$ is the Legendre symbol modulo $p$. The function $t(x)$ satisfies:

$$
t(x)= \begin{cases}1 & \text { if } \chi(x-1)=1, \chi(x)=-1, \chi(x+1)=-1 \\ 0 & \text { otherwise }\end{cases}
$$

Thus, the sum $S=\sum\limits_{x=1}^p t(x)$ counts the number of $x \in \mathbb{F}_p$ satisfying the conditions $\chi(x-$ $1)=1, \chi(x)=-1$, and $\chi(x+1)=-1.$

Expanding $t(x),$ we obtain
\begin{equation*}
    \footnotesize
   t(x) =\frac{1}{8} \left[ 1 + \chi(x\!-\!1) - \chi(x) - \chi(x\!+\!1) - \chi(x^2\!-\!x) - \chi(x^2\!-\!1) \!+ \!\chi(x^2\!+\!x) 
     +\chi(x^3\!-\!x) \right]. 
\end{equation*}
Note that the discriminants of \(x^2 - 1\), \(x^2 - x\), and \(x^2 + x\) are nonzero in \(\mathbb{F}_p\).
Combining the fact that \(\sum\limits_{x=1}^p \chi(x) = 0\), and applying  Lemma~{\rm{\ref{60}}},  we obtain that
$$S=
\frac{1}{8} \left( p+1+\sum\limits_{x=1}^p\chi(x^3-x)\right)
.$$
Applying Lemma~{\rm{\ref{50}}}, we have 
$S\geq \frac{1}{8} \left( \sqrt{p} - 1 \right)^2.$
This completes the proof of Corollary~{\rm{\ref{21}}}.
\qed
\section{\centering \textbf{Conclusion}}

\textbf{1.} Constructing an infinite family of irreducible polynomials over finite fields is an important topic in the area of finite fields (See \cite{Cheng}). We  achieve this in Corollary 2.8 (as detailed in the proof of Corollary 2.8).

\textbf{2.}
 The``Everywhere Eventual Stability Conjecture," as the name suggests, posits that there exist abundant rational maps defined over number fields or function fields that are eventually stable. However, this remains conjectural. In our work (Corollary {\rm \ref{37}} and Theorem {\rm \ref{88}}), we have found a class of inversely stable polynomials \(\phi(z)\) that hence naturally induce a family of eventually stable rational maps \(\left( \frac{1}{\phi(z)}, \infty \right)\). This construction thereby provides concrete evidence in support of the validity of the ``Everywhere Eventual Stability Conjecture."

\textbf{3.}
 From Corollary~\ref{37} and \cite{Rafe} Theorem~3.1, we obtain the following result.
 
Let \( S \) be a finite set of places of the rational number field \(\mathbb{Q}\) containing all archimedean places, let \(\infty \in \mathbb{P}^1(\mathbb{Q})\), and let \(\phi(z) = z^d + c \in \mathbb{Z}[z]\) be irreducible  \( d \geq 2 \). Let $\Phi(z)=\frac{1}{\phi(z)}.$ Suppose that either
\( d \) is odd, or \( d \) is even and \( c \) is not a square in \(\mathbb{Q}\).
Then, for every \(\gamma \in \mathbb{P}^1(\mathbb{Q})\) that is not preperiodic under \(\phi\), the intersection \(\mathcal{O}_{S, \gamma} \cap \mathcal{O}_{\Phi}^{-}(\infty)\) is finite, where \(\mathcal{O}_{\Phi}^{-}(\infty)\) denotes the backward orbit of \(\infty\) under \(\Phi\), and \(\mathcal{O}_{S, \gamma}\) is the ring of \( S \)-integers relative to \(\gamma\).

\textbf{4.} From Corollary~\ref{37}, Theorem {\rm\ref{018}} and \cite{Rafe} Proposition~2.2, we can obtain a coarse measure of the image size in the Arboreal Galois representation. (See \cite{Rafe} for details)

\end{document}